# Convergent Numerical Solutions for Unsteady Regular or Chaotic differential equations

Lun-Shin Yao
School for Engineering of Matter, Transport and Energy
Arizona State University, Tempe, AZ 85287

Abstract

Von Neumann established that discretized algebraic equations must be *consistent* with the differential equations, and must be *stable* in order to obtain convergent numerical solutions for the given differential equations. The "stability" is required to satisfactorily approximate a differential derivative by its discretized form, such as a finite-difference scheme, in order to compute in computers. His criterion is the necessary and sufficient condition only for steady or equilibrium problems. It is also a necessary condition, but not a sufficient condition for unsteady transient problems; additional care is required to ensure the accuracy of unsteady solutions.

**1. Introduction**

Systems of ordinary differential equations that exhibit chaotic responses have yet to be correctly integrated. So far no convergent computational results have ever been determined for chaotic differential equations, since the truncation errors introduced by discretized numerical methods are amplified for *unstable* computations. Numerical methods usually convert continuous differential equations to a set of algebraic equations to be solved by computers. Von Neumann established that discretized algebraic equations must be *consistent* with the differential equations, and must be *stable* in order to obtain convergent numerical solutions for the given differential equations. A typical property of chaotic differential equations is that they are *unstable*. It is not straightforward to check the consistence and stability of a numerical computation. In particular, it lacks a practical way to conveniently check the convergence of numerical results for *non-linear* differential equations that a linear stability analysis may not yield desirable and confident conclusions.

Parker and Chua [1] suggested a practical way of judging the accuracy of the numerical results from a non-linear dynamical system is to use two or more different methods to solve the same problem. If the two solutions agree then they can be assumed accurate. Viana [2] proposed to solve the same problem in two or more different machines to ensure the convergent results. Both approaches are testing to ensure that truncation errors will not overwhelm the correct solutions. The same propose can be achieved by solving the problem in one machine and one method, but two different integration time steps [3, 4]. All three ways are easy to apply, but the agreement of two computational results by either of these ways is only a *necessary* condition, and is not *sufficient*. Typical examples, well-known to all graduate students in thermal science, are unsteady heat conduction problems; even though, the heat equation is linear. They demonstrated the additional difficult of checking convergence for unsteady problems.

Without knowing it is not a sufficient condition, Lorenz [5] mistakenly concluded that his solution for his 1990 model was convergent initially for thirty years! This contradicts to the fact that the initial period of the Lorenz solution for his 1990 model is mixed with many *unstable and divergent* sections with some stable sections. One cannot claim that the mixture of errors in many unstable computations with some short time convergent computations is a correct solution, since the differential directives cannot be replaced by their computable discretized forms unstable periods. We will explain why the convenient ways to check convergence of unsteady computation is insufficient below and followed by numerical examples.

## 2. Mathematical Explanation

We will solve a set of, or a differential equation

$$\frac{du}{dt} = f(u;t), \tag{1}$$

whose exact solution is $u = u(t)$. Let's use $X_i$ (t) (i=1,2) denotes the computational results for two different methods, or two different machines, or two different integration time steps; $E_i$ (t) is the corresponding computational errors. Therefore,

$$X_i(t) = u(t) + E_i(t). \tag{2}$$

If the difference of two computational results is small, such as

$$|X_1 - X_2| = |E_1 - E_2| < \varepsilon, \tag{3}$$

where ε is a pre-assigned small number, it has a possibility that $E_1$ and $E_2$ are both small, and the computational results are convergent. On the other hand, (3) does not guarantee that both E's are small; it only states that the difference of two errors is small. Hence, (3) can only be a necessary condition. This is the mistake made by Lorenz [5].

## 3. Numerical Examples

Two examples will be given below to demonstrate the convergence of numerical solutions for differential equations.

A. The first one is a simple linear differential equation and we will construct stable computations to demonstrate that (3) is only a necessary condition.

The equation,

$$\frac{du}{dt} = -10u, \tag{4}$$

is used with the initial condition u(0)=1. The exact solution is

$$u = \exp(-10t). \tag{5}$$

The explicit finite-difference scheme is chosen for an unstable computation as

$$\frac{u^{n+1} - u^n}{\Delta t} = -10u^n,$$

or

$$u^{n+1} = (1 - 10\Delta t)u^n. \tag{6}$$

It is clear that (6) is unstable, if $\Delta t > 0.1$; the truncation errors is $O(\Delta t)$. It is well known that the numerical result will diverge for an unstable computation. We will show two computational results in comparison with the exact solution: one is for $\Delta t = 0.05$; the other $\Delta t = 0.06$. Two computational results agree completely initially. This is because that the computations of first step for two different time steps are identical. Since the truncation errors are $O(\Delta t)$, the results for the first few steps, when the computation time is about the same $O(\Delta t)$, cannot be accurate. The comparison presented in the Figure 1 demonstrates two computational results are not close to the exact solution even though they are fairly close to each other. This confirms the claim *that the small difference of two computational results can only be a necessary condition for the convergence of unsteady problems*. On the other hand, both computational results asymptotically converge to the exact solution, zero, for the steady problem. This is an example to show that von Neumann's consistent and stable conditions are necessary and sufficient for steady problems, but not sufficient for unsteady problems. For a *consistent* and *stable* computation, it still requires checking computed results by successively reducing time-step size until the difference is acceptably small; then, the convergence can be claimed [4] for an unsteady computation.

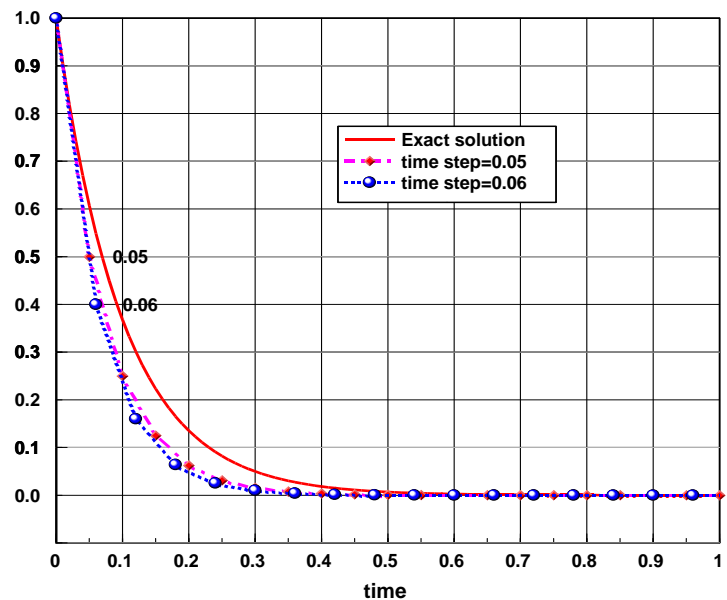

**Figure 1. Comparison of computational results with the exact solution. The number on the right of the plotting points indicates the computational time of the point.**

Another commonly known example, frequently taught in the first-year graduate course in heat transfer, is the heat equation. It is well known that a consistent and stable computation is sufficient to provide a convergent steady-state solution, but cannot guarantee a convergent transient solution. A convergent transient solution can only be obtained by successively reducing the integration time steps until the change of the computed transient results is acceptably small.

B. The second example is the Lorenz second model [4, 5, and 7]. The model is composed with three non-linear first-order differential equations.

$$\frac{dX}{dt} = -Y^2 - Z^2 - \frac{1}{4}(X-8),\qquad(7a)$$

$$\frac{dY}{dt} = XY - 4XZ - Y + 1,\qquad(7b)$$

$$\frac{dZ}{dt} = 4XY + XZ - Z.\qquad(7c)$$

The initial condition used below is (X=2, Y=1, Z=0). The error curve presented in Figure 2 is the difference of X(t) computed by the fifth-order Taylor-series method [7] with $10^{-6}$ time step by the Taylor-series method for time-steps $10^{-7}$, respectively. The conclusion is independent of the numerical methods used to integrate the equations (7). The details of comparison of various methods can be found in [4].

The error curve shown in Figure 2 differs obviously from any non-convergent error curves for any linear differential equations. The recorded difference of two computational results is too small when time is less than 30; so, we did not plot them. According to Lorenz's opinion [5], this shows that the numerical solutions are good for this short period of time; even thought, he agreed that numerical solutions for long time is not possible. It is worthy to point out that the time steps used in our computation is much smaller than what Lorenz used; so, our *good* results, according to Lorenz's criterion, can be extended to larger time. We will explain why this concept is wrong below.

The only available detailed error analysis for numerical solutions of non-linear differential equations, as we are aware, is for the famous Lorenz 1963 model [8]. It clearly demonstrated that two major amplification mechanisms exist for truncation errors, introduced by all numerical methods. The first is the explosive amplification mechanism, which can instantly amplify the truncation tremendously when the trajectory penetrates the separatrix by violating the differential equations. Since the Lorenz 1990 model does not have an attractor, the explosive amplification does not occur, confirmed by our numerical computations [4]. We will not further discuss it here; the interested readers can read [8, 9].

The second mechanism is the exponential amplification of errors, which is also found in the numerical solution of linear differential equations as explaining in the first example. An unstable computation for unsteady linear differential equations can result two kinds of behaviors uniformly in time: exponential growth of errors, or exponential growth of the amplitude of oscillatory solutions. The crucial difference between non-linear differential equations and linear ones is the exponential error amplification for non-linear differential equations is not uniform in time, see Figure 2. The growth of truncation errors occurs in "irregular valleys." This suggests the existence of certain dynamic structures in the phase space. This agrees with the *exponential amplification* of errors described in [8, 9]. When two trajectories move along the direction of a stable manifold, the distance between them shrinks; in other words, errors are reduced; when trajectories move along an unstable manifold, errors are amplified. The combined consequence is, however, the exponential growth of truncation errors in time as shown in Figure 2.

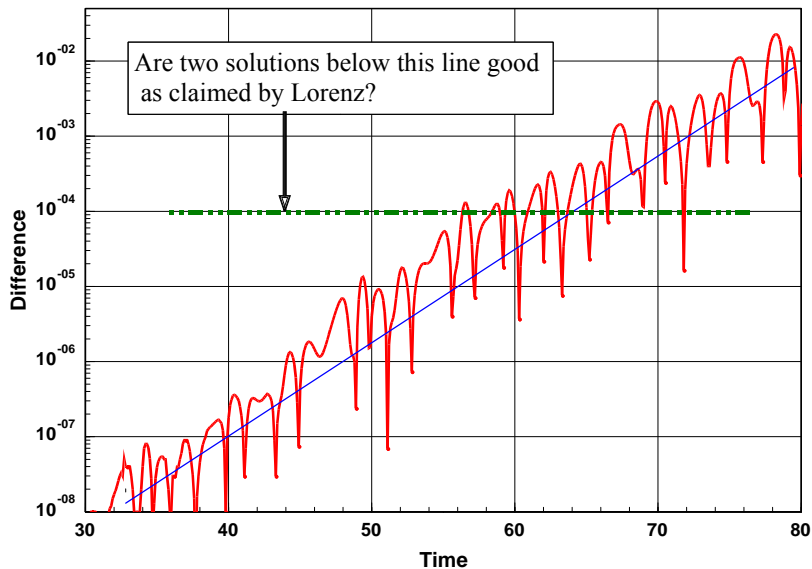

**Figure 2. The difference of numerical results of time steps** $10^{-6}$ **and** $10^{-7}$.

It should be emphasized here that the error amplification is due to the unstable computation locally, which violates von Neumann's convergence criterion. For linear differential equations, it will lead to divergent solutions; one would not expect it could provide correct solutions for non-linear differential equations. If it were so, it would mean that it were easier to numerical integrated non-linear differential equations, since no check of convergence would be needed. Then, it was always legitimate to replace a derivative by a finite-difference counterpart without worrying they may not even be approximations! This is exactly what has happened in solving chaos or turbulence numerically now.

It is also worthy to mention that it has been demonstrated in [8, 9] that a small difference of two computations does not imply either one is close to the correct solution for unsteady problems. This has been experienced many times in the history of numerically solving both linear and non-linear differential equations of unsteady problems, but has been overlooked in solving chaos or turbulence numerically.

This difficulty associated with unstable computation is the property of non-linear differential equations, and cannot be remedied by adjusting numerical methods, see [8, 9]. Since the truncation errors are not controllable and occur randomly, the numerical computational chaos results, or turbulence is also random in nature; irrespectively, the associate boundary conditions are either independent of time, or depend on time regularly. Consequently, an unstable numerical result is the random amplification of truncation errors, induced by numerical processes, and has no physical meaning.

We do not believe that our paper can reverse the avalanche of treating numerical errors as numerical solutions of differential equations, but hope someone, in the near future, may take a little effort to honestly compare computational results with carefully carry-out measurements. It is time to reconsider the activities of continuously producing numerical errors in large amount without any justification. Fundamental principles in science should always be respected before one can prove otherwise.

**4. Comments, discussions, and open problems from some experts in the field:**

I will first outline some well-known basic principles in numerical mathematics, which will help to explain the following questions.

In calculus, we know

$$\frac{\partial u}{\partial t} = \lim_{\Delta t \to 0} \frac{u(t+\Delta t) - u(t)}{\Delta t}. \tag{8}$$

Since computers are digital computational devices, it is necessary to *discretize* a derivative in order to calculate it in computers. I will use a finite-difference scheme as an example below; the principle can be equally applied to all discretized numerical methods without exception. A derivative can usually be replace by a finite-difference form,

$$\frac{\partial u}{\partial t} = \left[\frac{u(t+\Delta t) - u(t)}{\Delta t}\right] + \text{TE} \ (\textit{truncation errors}), \tag{9}$$

where TE represents truncation errors and is of $O(\Delta t)$ in the above example. In the limit of $\Delta t$ approaches zero, (9) agrees with (8). The TE always exists in any discretized process. Algebraic equations are resulted after all derivative terms being replaced by their difference forms; for example, terms in the square bracket in (9). If the resulting algebraic equations are *stable*, TE will be exponentially decayed and the term inside the square bracket is a good approximate of the derivative, since (8) will be asymptotically satisfied. On the other hand, if the resulting equations are *unstable,* the TE will be exponentially growing, and (8) is violated, or it implies that

$$\left[\frac{u(t+\Delta t) - u(t)}{\Delta t}\right] \text{ does not approach } \frac{\partial u}{\partial t}. \tag{D3}$$

The resulting equations after discretization have nothing to do with the original differential equations, and certainly are not an approximation of the original differential equations. Thus, the solutions of such algebraic equations are unrelated to the original differential equations. For linear differential equations, an unstable computation usually results exponentially divergent results, or exponentially divergent oscillatory results, a clear sign of a failure computation. This is why von Neumann put forward that the "stability" is the necessary and sufficient condition for a convergent solution for *steady-state* problems.

For *unsteady* problems, it is necessary to use even smaller $\Delta t$ than required the one by Courant-Fredrich-Levy (CFL) condition in order to get an accurate transient solution as demonstrated in the first example of the current paper. There are many examples can be found in heat conduction problems in undergraduate heat transfer textbooks. Consequently, stability is only a necessary condition to obtain an accurate transient solution.

I will answer and discuss the following questions.

1. The author comments the stability of the numerical computations for nonlinear ordinary differential equations. More precisely, the author comments that the agreement of two computational/numerical results is only necessary not sufficient condition for the accuracy of the computational/numerical results. Really, the accuracy of the computational/numerical results is very important scientific problem. Unfortunately, the author does not propose a way for the solution of this problem.

Ans: The answer to this question is very simple for linear or non-linear differential equations, if the resulting algebraic equations are *stable*; continuously reducing the integration time step, $\Delta t$ until the change of the transient solutions for two different time steps becomes acceptably small. Then, both results can be considered as an accurate transient solution. If one carries out the computation of the first example in the paper by further reducing the integration time step, an accurate transient solution can be readily found, and agree with the exact solution to any degree as one wish. There are many other examples of heat conduction problems can be found in undergraduate heat transfer texts.

For non-linear chaos differential equations, *when the governing parameter is larger than its critical value*, the situation becomes much more complex. Two Lorenz's models have been discussed in [4, 8] and his second model is also used as the second example of this paper. The Lorenz's first model was analyzed in detail and reported in [8]. It shows that the truncation (numerical) errors are amplified exponentially in the unstable region (manifold); are reduced exponentially in the stable region (manifold). According to the basic principle outlined above, it is clear that the numerical results in the unstable regions cannot be considered as an approximate solution to the original differential equations. Similar conclusion can be made for the Lorenz's second model and reported in [4].

In addition, we have identified that the Lorenz's first model contains local separatrix, not in his second model. The truncation errors can be amplified explosively when the trajectory penetrates the virtual separatrix, which violates the differential equations. The existence of a virtual separatrix is a consequence of singular points of a non-hyperbolic system of differential equations, which is not *shadowable* [8, 9]. A commonly cited computational example in chaos involving two solutions of slightly different initial conditions that remain "close" for some time interval and then diverge abruptly when one penetrates the virtual separatrix, violating the differential equations and the other does not. Before it was pointed out in [8, 9], that this phenomenon is actually due to the explosive amplification of numerical errors, and violation of the differential equations as described above, this behavior was often believed to be *a typical characteristic of chaos*, and frequently used as the evidence that a computation is chaos. Similar computational results can occur for two different integration time steps; many would consider such results as acceptable since it is a "twin brother" of sensitivity to initial conditions and a typical characteristic of chaos. This mistake deserves clarification and has been explained thoroughly in [8].

Additional examples for other non-linear chaos differential equations are cited in [3].

The central issue discussed in [9] is that no chaos solution exists for differential equations, since all computations are unstable. This is the obvious consequence that all discretized numerical methods have truncation errors and are incapable to solve chaos differential equations.

Since non-linear differential equations can have *multiple solutions* when the value of its parameter is larger than its critical value, once the trajectory went in the unstable region, the amplified truncation errors alters its initial condition equivalently for the next stable region. This can lead to next stable solution from the one that one was originally trying to get as discussed in [D1]. This is a brand new topic, which has never been studied, yet.

2. I tried to re-evaluate the revised version of the manuscript in goodwill. However, in his revised manuscript the author insists in keeping his claim that "Lorenz [5] mistakenly concluded that his solution for his 1990 model was convergent initially for thirty years!" (in Introduction). For me this was the one and only major issue I had, as reading Lorenz's answer to Yao and Hughes published in Tellus (2008), 60A, 806--807, I think that this claim is both strong and wrong.

Ans: The answer of this question is simple and clear by judging it from the basic principle outlined above. Lorenz's numerical results of his second model clearly show the trajectory went through stable and unstable regions alternatively, see Figure 2 of the paper and [4]. The difference of the results for two different time steps increases when the trajectory moves in the unstable region, where (8) is violated; the difference decreases when the trajectory moves in the stable region. Can one claim such a computational result is good when most parts of it violate (8), a definition in calculus?

Many authors followed Lorenz's step [5] and claimed mistakenly that their solution is good for a short initial period. Form the above discussion; it is clear that the initial good period will be zero, if the initial point is selected in the unstable region. If the initial condition is located in the stable region, the computed result will be *good* until it moves out of the stable region and gets into the unstable region.

Since without numerical instability, there is no chaos; so they can only claim their results is a good regular solution, not chaos for the initial period. I want to emphasize again that all numerical chaos are amplified numerical errors.

3. The topic of this paper is rather interesting. However, the extension and scope of the paper seems to be very brief for what such a topic may require. The references seem to be devoted to self-citations, when there is a huge amount of scientific literature in the field.

Ans: Detailed analyses are available in [3, 4, 8, 9] and additional references are listed below, which are all my own work. I do not know the existence of any written paper in line with my conclusion. If there are some, I welcome readers to reveal them to me. The entire huge amount of scientific literature devotes to treat amplified numerical errors as computed chaos with very few exceptions cited in the additional references listed below. This is exact the reason that I tried to push this paper forward. Fortunately, science is not a democratic system that majority wins; in fact, truth always prevails in science. I hope the basic principles, outline above, can help readers to analyze the situation rationally by themselves, not just follow the argument of established authorities.

4. The major part of the discussion seems to be devoted to demonstrate "that the small difference of two computational results can only be a necessary condition for the convergence of unsteady problems". All the discussion seems to emerge from Eq.3, which from another side (may be a too simplistic one), seems to be self-evident: "getting two bad results rather similar does not guarantee that the conclusion is alright". The focus of the discussion could be set in dealing with the definition of good or bad.

Ans: I do not really understand this question, but I believe the above discussion has already answered this question.

5. It may be of interest to cite some of the large amount of references in the literature devoted to the topic of solving a problem by computing many similar initial conditions/parameters, and analysing later the results. This aims to assign given probabilities of success to every computed solution (as for instance in the field of meteorology).

Ans: Please read the answer of the second question, and note that the large amount literatures ignored the basic principles discussed above.

I would like to point out that a slightly different initial condition can result a completely different long-time solution for non-linear differential equations when the governing parameter is above its critical value; so there are multiple solutions exists even for non-turbulent flows. I list some references [9-12] below for interesting readers to explore. More references can be found in them. They are also all my works, since my works are the only theoretical studies exist, beside few experiments in fluid mechanics.

It is worthy to mention, from our experience, that the computation of the Navier-Stokes equations fails to converge if the Reynolds number is too much larger than the critical Reynolds number. This is the reason why a direct numerical simulation cannot be use to study flow transition; so, it is also unhelpful for turbulence.

A well-known excise to stabilize an unstable computation is to use upwind difference scheme; however, the added *numerical viscosity* associated with the upwind difference scheme may overwhelm the actual viscosity to invalidate the computational results. One can see why weather forecast is so unreliable, it is not because that the Navier-Stokes equations are incorrect; it is due to lack of a method to solve them correctly.

6. Along the paper, the concepts of chaoticity and stability seem to be mixed somehow, and a clear separation of both definitions should be appreciated.

Ans: The rigorous mathematical definition of chaos can be found in [2, 8].

Traditionally, the linear stability analysis in fluid dynamics is to study the growth or the decay of a very small perturbation quantity added to the steady-state base flows. For a laminar flow, the differential equations are stable, but the difference scheme can be unstable. This leads to CFL condition.

The extension of such analyses for time-dependent base flows is complex and unsuccessful. The linear-stability analysis of the algebraic equations resulted from the discretization of differential equations is very similar to the stability analysis in fluid dynamics. An easy alternative is to solve differential equations with two different time steps, and hope the difference is acceptably small as discussed in this paper.

For chaos differential equations, it becomes very complex, since the differential equations are themselves unstable. There is no way to design a stable numerical method for unstable differential equations. As noted in the Question 5, there is a well-known way to stabilize an unstable computation by using use upwind difference scheme; however, the added *numerical viscosity* associated with the upwind difference scheme may overwhelm the actual viscosity to invalidate the computational results.

I have repeatedly tried to explain that numerical chaos is simply amplified numerical errors. Numerical chaos and numerical instabilities are different titles, but the same object.

7. The last paragraph on page 6 deals with the numerical computational chaos results. I agree with the fact that any numerical scheme will diverge from a true orbit beyond certain timescales for given problems. But some discussion about how these timescales may vary depending on the nature of the orbit, and even may be very long even when the orbit is chaotic is of interest. A discussion about the shadowing and predictability topics should be also appreciated.

Ans: This answer of this question is a part of the Question 2 and copy below for the convenience of readers.

Many authors followed Lorenz's step [5] and claimed mistakenly that their solution is good for a short initial period. Form the above discussion; it is clear that the initial good period will be zero, if the initial point is located in the unstable region. If the initial condition is selected in the stable region, the computed result would be *good* until it moves out of the stable region and gets into the unstable region.

The concept of shadowing is briefly reviewed in [9]. It was originally invented to save chaos theories for the hyperbolic systems. Since the relation between chaos theory and differential equations has not been established (Smale's $14^{th}$ problem) and topological transitivity cannot be proved and is likely invalid, shadowing is not a useful concept for differential equations.

8. There is a vast amount of literature devoted to the numerical methods and how they deal with chaos. A brief panorama of the field should be of interest. After this, it may happen the statement "but hope someone, in the near future, may take a little effort to honestly compare computational results with carefully carry-out measurements" should be revised by the author.

Ans: It is a fact there is a large amount of literature devoted to the numerical methods and chaos, since both are the main streams of research in their areas and have been heavily funded by government agents; in particular, in the US. How could so many smart researchers all made the same mistake of violating the basic principles of numerical methods in solving differential equations, put forward by von Neumann?

Maybe, this is due to the limitation of human brains to discover brand new idea; instead, we look established authorities for guidance. However, please read my answer of the question 3 above.